\documentclass{gazette}

\input{psfig.sty}

\usepackage{amssymb}

\def\Z{\mathbb{Z}}
\def\round#1{\left\lfloor#1\right\rceil}

\begin{document}

\title{Graphical patterns in quadratic residues}

\author{I. Jim\'enez Calvo}
\address{C.S.I.C., C/Serrano 144, 28006--Madrid, Spain.}
\email{ijcalvo@terra.es}

\begin{abstract}
Plots of quadratic residues display some visual features that are analyzed mathematically in this paper.
The graphical patterns of quadratic residues modulo $m$ depend only on the residues of $m$ modulo the 
lowest positive integers.
\end{abstract}

\maketitle

\section{Introduction}
The quadratic residue (QR) of an integer $x$ modulo $m$ is the remainder $a$ 
of the division of $x^2$ by $m$. It is expressed by the congruence,
$$
x^2 \equiv a \pmod m,
$$
where the residue class $a$ can be represented by an integer from 0 to $m-1$. 
Figure 1 shows a typical plot of QRs where a rather chaotic distribution can be appreciated. 
As $x$ and $m-x$ have the same QR, then only those QRs corresponding to $0 \le x < m/2$ are represented.
It is not surprising that QRs modulo a big integer of unknown factorization are used in cryptography 
and random sequence generators \cite{BBS}. Nevertheless, they always show a remarkable 
presence of points that seem to draw multiple ``acute parabolas'' that spread along the graphic.
A closer look confirms that the acute parabolas have a quite regular distribution, for values of $x$ 
near simple fractions of the modulus in a quantity related to the denominator of the fraction, 
so that a lower number of 
best defined acute parabolas correspond to a simpler fraction .  Not only 
QRs, but also $n$-power residues spread as integer points of sets of $n$-degree 
polynomials for $x$ values near a fraction of the modulus as it is stated in \cite{JS}. In that paper, 
it was shown that some of those ``acute parabolas'' can be used to obtain low quadratic residues 
used in the factorization of large integers. In fact, those parabolas correspond to the 
Montgomery polynomials 
of the MPQS ({\it Multipolynomial Quadratic Sieve}) \cite{Pomerance,Silverman}  and 
HMPQS ({\it Hypercube Multipolynomial Quadratic Sieve}) \cite{Peralta}  factorization algorithms.

A detailed treatment of the presence and regular disposition of the ``acute parabolas'' is 
done in Section 2. Nevertheless, such features can be explained in a rather informal way with 
a fairly straightforward reasoning. Let be an integer $x$ and a modulus $m$ and consider the 
difference between the square of $x$ and its nearby integers, i.e.  
$\Delta=(x+i)^2-x^2=2ix+i^2$, when $x$ is the closest integer to a fraction $a/b$ of the modulus. 
We consider small $i$ values and neglect the term $i^2$. A second approximation is done replacing 
$x$ by $(a/b)m$, so $\Delta \approx  (2ia/b) m$.
If $b$ is odd, the residue of the fraction $(2ia/b) m$ modulo $m$ is another fraction 
$ (a'/b) m$ with $a' \equiv 2ia \pmod b$ because $2ia=a'+kb$ for some integer $k$. 
When $b$ is even, we can note that $\Delta \approx (ia/b') m$ with $b'=b/2$ 
and its residue modulo $m$ is the fraction $ (a'/b') m$ with $a' \equiv ia \pmod {b'}$.
Therefore, the QRs for values near $x$ are spaced approximately by multiples of $m/b'$. 
The term $i^2$, which is not considered in this approximation, accounts for building the parabolas at 
regular intervals.

Functions involving QRs have been used to generate computer graphics. The function 
$f(x,y) = x^2+y^2 \pmod m$ shows concentric circles with centres at coordinates that are made up of simple fractions of the modulus. They are actually, truncated paraboloids present 
in this bi-dimensional case (see fig.~2). Another type of patterns can be observed sometimes in QRs 
plots besides the ``acute parabolas''.  For some moduli, it may happen that the outlines of broader 
diffused ``ghost parabolas'' materialize, as can be observed in fig.~1 and fig.~3. 
These kind of patterns are analyzed in Section 3. To conclude, we can state that the graphical 
patterns of QRs mod $m$ plots depend only on the residues of $m$ modulo the lower positive integers.

\begin{figure}
\centerline{\psfig{figure=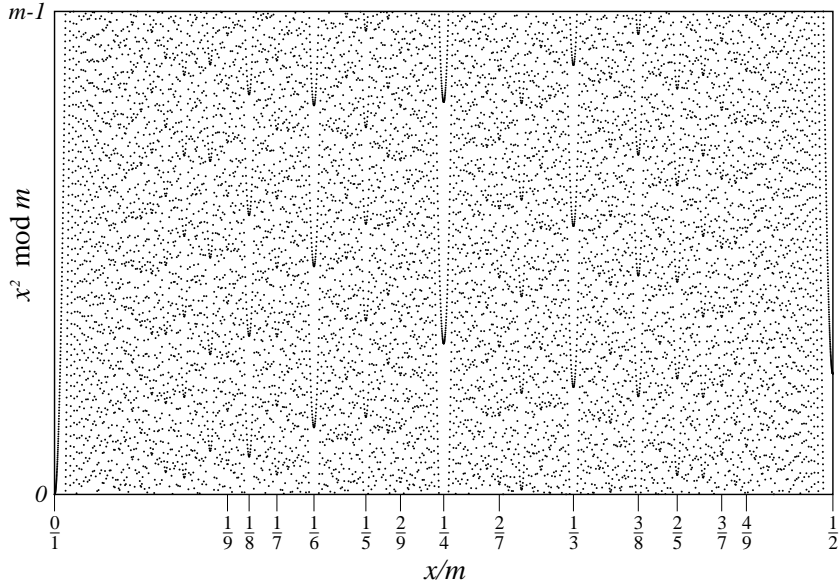,height=8cm}}
\caption{Quadratic residues mod 20171.}
\end{figure}

\begin{figure}
\centerline{\psfig{figure={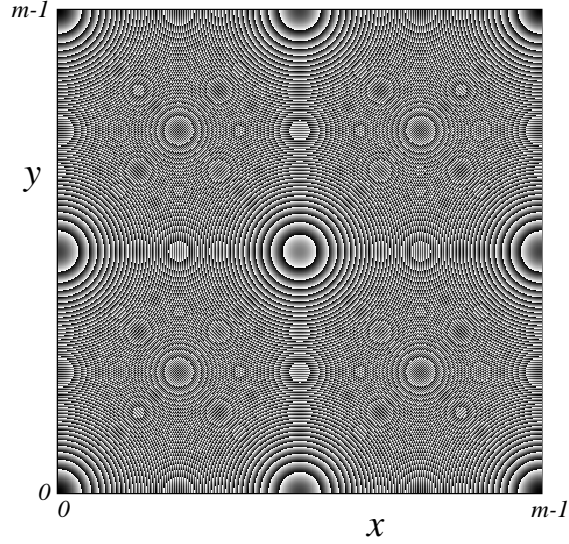},height=8cm}}
\caption{Function $f(x,y)=(x^2+y^2)$ mod 415.}
\end{figure}

\section{The geometry of the ``acute parabolas''}
The congruence relation is a relation of equivalence between integers. Nevertheless, a similar relation 
between rationals can be established:
\begin{definition}
Two rationals, $s$ and $t$, are Q-congruent modulo an integer $m$ when $(s-t)/m$ is an integer 
and then we can write
$$
   s \equiv_Q t \pmod m\quad \mbox{and}\quad s = t + km,
$$
for some integer $k$.
\end{definition}
$Q$-congruence relation is reflexive, symmetric and transitive as well as the congruence relation 
between integers. Moreover, two congruent integers are also $Q$-congruent. Therefore, the notation 
``$\equiv$'' between rationals must be read as ``$\equiv_Q$'' and when the relation 
is between integers, it can be read in both ways.

Let $m$ be a modulus and $a/b$ an irreducible fraction. We define the integer $x_0$ and its 
quadratic residue $r_0$ as  
\begin{equation}
   x_0=\round{ \frac a b m },\quad x_0^2 \equiv r_0 \pmod m \; ,
\end{equation}
where $\round{x}=\lfloor x+1/2\rfloor$, that is, the nearest integer to $x$.

\begin{proposition}
$r_0$ is the nearest integer to the rational $\beta (m/b^2)$ where \\
$$
\beta \equiv a^2m - 2a\alpha \pmod {b^2} \quad {with} \quad \alpha \equiv am \pmod b, \:|\alpha| \leq b/2.
$$
\end{proposition}
\begin{proof}
We note that $x_0=(am-\alpha)/ b$ from the definition of the integer $\alpha$ and then, 
$$
  x_0^2=\frac {a^2m^2-2a \alpha m} {b^2} +  \frac {\alpha^2} {b^2} =
\left({\frac {a^2m-2a\alpha-\beta } {b^2}}\right)m+ \beta \frac m {b^2}+
 \frac {\alpha^2} {b^2} \; .
$$
With the definition of $\beta$, the term $(a^2m-2a\alpha-\beta)/b^2$ becomes an integer. Taking 
into account that $|\alpha| \leq b/2$, we can then write
\begin{equation}\label{r_0}
  r_0 \equiv \beta \frac m {b^2} +\frac {\alpha^2} {b^2} \pmod m,
  \quad 0 \leq \frac {\alpha^2} {b^2}  \leq 1/4 \;.
\end{equation}
\end{proof}

We now calculate the QRs of integers near $x_0$. From the considerations explained in 
the introduction, it becomes natural to consider lattice points defined by
\begin{equation}\label{lattice}
x=x_0+i+jb', \quad {\rm with} \quad b'=\left\{
 \begin{array}{lc}
    b, & {\rm if\:{\it b}\: odd}\\
    b/2, & {\rm if\:{\it b}\: even}
 \end{array} \right. ,\quad |i| \le b'/2, \quad j \in \mathbb{Z}_m.
\end{equation}
Expanding the square of the trinomial, we have
$$
x^2 = (x_0+i)^2+j^2{b'}^2+2(x_0+i)jb'.
$$
Taking into account that $x_0=(am-\alpha)/ b$ and with $r$ and $r_i$ as the quadratic 
residues of $x$ and $x_0+i$ respectively, we can write
\begin{equation}\label{r}
r \equiv {b'}^2j^2+2\left(b'i- \frac \alpha c\right)j+r_i \pmod m \;,\quad
c=\frac b {b'} = \left\{
\begin{array}{ll}
    1, & {\rm if\:{\it b}\: odd.}\\
    2, & {\rm if\:{\it b}\: even.}
 \end{array} \right. 
\end{equation}

\begin{proposition}
The quadratic residues of integers near $(a/b) m$ split as integer points of 
$b$ parabolas if $b$ is odd, or $b/2$ parabolas if $b$ is even, with their
vertices placed at $x_v=(a/b) m$ and equidistant in between.
\end{proposition}
\begin{proof}
It is clear that (\ref{r}) is a parabola under variable $j$ for each one of the $b'$ values 
of $i$ in (\ref{lattice}). Those parabolas will be named as $P_i(j)$. 
The coordinates at the vertex are given when the derivative of (\ref{r}) with respect to $j$ becomes null. 
Since $r_i$ is fixed for each $i$ value they are given by 
$$
j_v= \frac \alpha {c{b'}^2}-\frac i {b'}, \quad
P_i(j_v)=r_i-\left(i-\frac \alpha b\right)^2  \;.
$$
Substituting the value of $j_v$ in (\ref{lattice}), and taking into account that $x_0=(am-\alpha)/b$, 
the coordinate at the vertex is $x_v=(a/b) m$ which is independent from $i$ value.
Then, all parabolas related to the fraction $a/b$ have vertices at the same $x$ value.

Since $(x_0+i)^2=x_0^2+2ix_0+ i^2$, then  
$r_i \equiv r_0 + 2i(am - \alpha)/b + i^2  \pmod m$.
Replacing $r_i$ and simplifying we obtain,
$$
P_i(j_v) \equiv r_0+2ia \frac m b- \frac {\alpha^2} {b^2} \pmod m\;.
$$
Taking the value of $r_0$ from (\ref{r_0}), we obtain, 
$$
P_i(j_v) \equiv \beta \frac m {b^2} + 2ai \frac m b \pmod m\;.
$$
Let be 
$$
a' \equiv \frac 2 c ia \pmod {b'}, 
$$
with $c$ as defined in (\ref{r}).  Then
\begin{equation}\label{minimos}
P_i(j_v) \equiv \beta \frac m {b^2} + a'\frac m {b'} \pmod m\;.
\end{equation}
As much as $i$ runs over $\Z_{b'}$ and $\gcd(a,b')=1$, $a'$ takes all values in $\Z_{b'}$.
From the above equation, it follows that vertices of $P_i(j)$ 
lie exactly at multiples of $m/b^2$ and are separated by exact gaps of $m/b'$. 
This concludes the proof. 
\end{proof}
 
\section{Patterns}
Propositions 1 and 2 give account of the number and position of the ``acute parabolas''. 
Near $x=(a/b) m$, $b$ parabolas (or $b/2$ if $b$ is even) appear and all are equidistant in 
$\Z_m$. Figure 2 shows a representation in gray levels of the function $f(x,y)=(x^2+y^2) \pmod m$. 
It is remarkable the symmetry of the graphic because it displays all the range of variables in $Z_m$. 
There is only one paraboloid at each node 
related to fractions 0/1, 1/2 and 1/1. That corresponds to centre, corner and the midpoint of the sides 
of fig.~2. 
Such paraboloids have a sharp profile form and they are surrounded by rings due to truncations.
Two paraboloids are present at each node related to 1/4 and 3/4 producing more fuzzy patterns. 
Other patterns alike can be observed for more complex fractions. 

At first glance, when plots like fig.~1 are observed, the presence of ``acute parabolas'' and 
their positions can be recognized. 
Since $b= b'c$, equation (\ref{minimos}) can be rewritten as
$$
P_i(j_v) \equiv \frac \beta {cb} \frac m {b'} + a'\frac m {b'} \pmod m.
$$
Let $\beta' \equiv \beta \pmod {cb}$, then
$$
P_i(j_v)  \equiv \frac {\beta' +sbc} {bc} \frac m {b'} + a' \frac m {b'} \pmod m,
$$
for some integer $s$. 

Taking $k=s+a'$, we have
\begin{equation}\label{vertices}
P_i(j_v)  \equiv \beta' \frac m {b^2} + k \frac m {b'} \pmod m,\quad  k \in \Z_{b'}.
\end{equation}

If the set of parabolas $P_i$ corresponding to the fraction $a/b$ is considered as a whole, without paying 
attention to the position of any single parabola in the set, the above equation shows that their 
position is defined only by the value of $\beta'$. Proposition 1 shows that $\beta$ 
(and also $\beta'$) depends only on $m$ for each fraction $a/b$. Since the sets of parabolas 
corresponding to the simplest fractions are those that are relevant to the eye, 
it follows that the general appearance 
of a plot only depends on  the values 
of $m\; {\rm mod}\; {cb}$ for the lowest values of $b$, that is for $cb=4,3,8,5,12, \cdots$. 

Let us consider the integer
$$
\Lambda(n)= 2 \,{\rm LCM}(2,3,4\cdots n).
$$
For example, let $\Lambda=\Lambda(9)=2^4\cdot3^2\cdot5\cdot7=5040$. Plots of QRs modulo integers that are congruent 
modulo $\Lambda$ will have the parabolas corresponding to fractions with denominator equal or 
less than 9  in the same position. Moreover, since $cb$ divides $\Lambda$ 
for $b=10,12,14,15,16,18 \cdots$, the parabolas related to fractions with these denominators will 
also have the same position. 
As an example, Figures 3(a) and 3(b) have the same appearance, but for their density, 
because the modules are $4\Lambda+19$ and $5\Lambda+19$ respectively and then 
$\beta\; {\rm mod}\; {cb}$ do not change when $b\le 9$ or $b=10,12,14,15,16,18 \cdots$.
Furthermore, we can prove that the vertices are rational points of a bundle of lines.
\begin{proposition}
Let the integers $\Lambda$, $m$ and $s$ be such that $s \equiv m \pmod \Lambda$ and let ${\mathcal B}$ be 
the set of positive integers $b$ such that $b|\Lambda$ if $b$ is odd or $2b | \Lambda$ if $b$ is even. 
The vertices of the parabolas defined in Proposition 2 associated to fractions $a/b$ with $b \in {\mathcal B}$ 
are rational points of the bundle of lines of the form
$$
   Y \equiv 2nX-sX^2 \pmod 1,\quad n=0,\pm1,\pm2, \cdots,\quad X=\frac x m,\;Y=\frac {x^2 \pmod m} m.
$$
\end{proposition}
\begin{proof}
The relative position of the vertices depends on the value of $\beta$ mod  $cb$ as it was justified above.
When $b$ is odd, $c=1$ and the expression of $\beta$  in Proposition 1 can reduce to 
$\beta \equiv -a^2m \pmod b$. If $b$ is even, $c=2$. We put $\alpha = am -kb$, for an integer $k$ 
such as $|\alpha| \le b/2$. 
Then, $\beta \equiv a^2m-2a(am-kb) \pmod {b^2}$, that reduced to mod $2b$ is also 
$\beta \equiv -a^2m \pmod {2b}$. We conclude that $\beta \equiv -a^2m \pmod {cb}$ for any parity of $b$. 
Note that, from the definition of $\Lambda$ 
and $s$, we have that when $b\in{\mathcal B}$, $s \equiv m \pmod {cb}$ and then 
$\beta \equiv -sa^2 \pmod {cb}$.
From Proposition 2 and equation (\ref{vertices}), it follows that the vertices of the parabolas are at 
$$
X_v=x_v/m=a/b,\quad Y_v=P_i(j_v)/m \equiv (\beta\;  {\rm mod}\; cb)\frac 1 {b^2} + k \frac 1 {b'} \equiv
-s \frac {a^2} {b^2} + \frac {ck} b \pmod 1.
$$
Note that $ck \equiv 2na \pmod b$ can always be solvable in $n$ for any value of $k\in \Z_{b'}$ 
and for any parity of $b$. Note also that $n$ can take any value with the appropriate value of 
$k\in \Z_{b'}$. Then, we have that 
$Y_v=-s \frac {a^2} {b^2} + 2n \frac a b = -sX_v^2 + 2nX_v$ for some $n$. 
\end{proof}

When $s=0$ the bundle of quadratics degenerates in straight lines( fig.~3(c)) that are outlined due to 
the accumulation of points near the vertices. As $s$ has a higher absolute value, 
the straight lines bend into parabolas as can be seen in fig.~3 (d) and 3(e).
When $s$ has absolute value as high as in fig~3(f), the lines of the bundle have such a high slope 
that they do not allow to notice the ``ghost parabolas''.

To conclude, we have shown that Proposition 1, 2 and 3 explain the main visual features of QR plots 
and that those features are related solely with the residues of $m$ modulo the lower positive integers.

\begin{figure}
\centerline{\psfig{figure={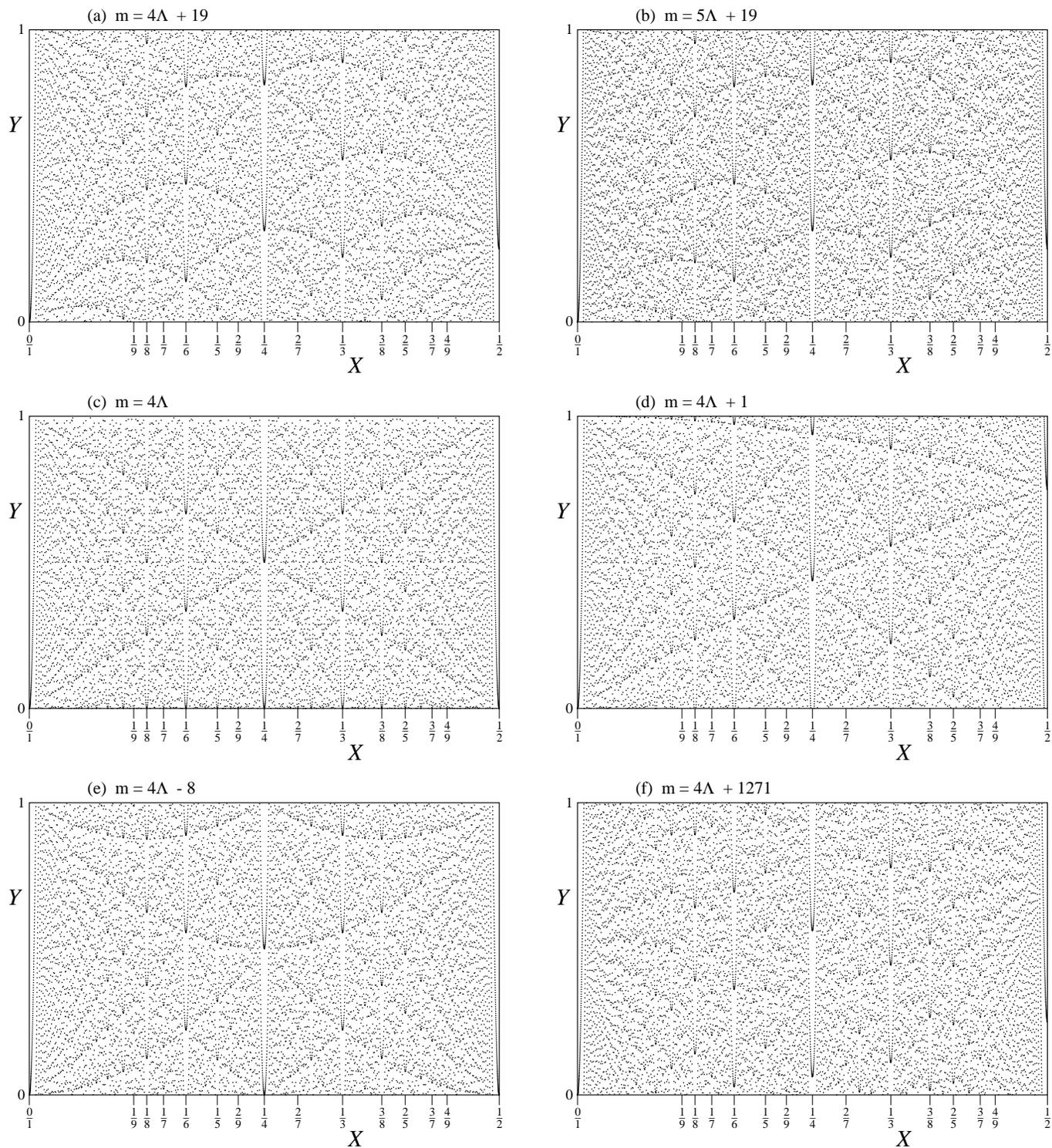},height=20cm}}
\caption{Plots of quadratic residues, $X=x/m$, $Y=(x^2\; {\rm mod}\; m)/m$, $\Lambda=5040$.} 
\end{figure}

\end{document}